\documentclass[11pt,a4paper]{article}

\usepackage[T1]{fontenc}
\usepackage[utf8]{inputenc}
\usepackage{lmodern}
\usepackage{microtype}
\usepackage{amsmath,amssymb,amsthm,mathtools}
\usepackage{thm-restate}
\usepackage[margin=3cm]{geometry}
\usepackage{enumitem}
\usepackage[hidelinks,hypertexnames=false]{hyperref}
\usepackage[nameinlink,capitalise,noabbrev]{cleveref}
\usepackage{bbm}

\setlist{nosep}
\allowdisplaybreaks
\usepackage[style=numeric,sorting=none]{biblatex}
\newtheoremstyle{mytheorem}
  {3pt}
  {3pt}
  {\itshape}
  {}
  {\bfseries}
  {.}
  { }
  {}

\newtheoremstyle{mydef}
  {3pt}
  {3pt}
  {\normalfont}
  {}
  {\bfseries}
  {.}
  { }
  {}

\theoremstyle{mytheorem}

\newtheorem{Th}{Theorem}[section]
\newtheorem{Lemma}[Th]{Lemma}
\newtheorem{Prop}[Th]{Proposition}

\theoremstyle{mydef}

\newtheorem{Def}[Th]{Definition}
\newtheorem{Conj}[Th]{Conjecture}
\newtheorem{Rem}[Th]{Remark}
\newtheorem{Prob}[Th]{Problem}

\title{\textbf{When Strategies Do Not Determine Play: Infinite Hat-Guessing Games on Subsets of the Reals}}
\author{Dániel Hegedűs\thanks{Leipzig University; Max Planck Institute for Mathematics in the Sciences, Leipzig; Eötvös Loránd University, Budapest. Email: \texttt{daniel.hegedus@mis.mpg.de
}}
\and
Benedek Kovács\thanks{HUN-REN Alfréd Rényi Institute of Mathematics, Budapest; Eötvös Loránd University, Budapest. Email: \texttt{kovacsbenedek@renyi.hu}}}
\date{}

\makeatletter
\renewcommand{\@fnsymbol}[1]{\@alph{#1}}
\makeatother

\begin{document}

\maketitle

\begin{abstract}
We study infinite hat-guessing games indexed by linearly ordered sets $S\subseteq\mathbb R$, where players guess sequentially using information determined by the order. While well-ordered standing sets admit the naive, recursive definition of play, we found a surprising example showing that this approach fails for non-well-ordered ones, where local decision rules do not determine a unique global outcome. By requiring uniqueness of play, we aim to provide a natural framework for hat-guessing games on arbitrary standing sets.

For the classical visibility model, it was previously known due to Glivický that at most one incorrect guess can be guaranteed when $S$ is well-ordered. We prove the same result for backwards well-ordered sets, even in the case when $S$ has no first element.

In the full-visibility version of the game, Glivický proved that at most one incorrect guess can be guaranteed if the standing set has a first element. We extend this result to arbitrary nonempty standing sets.

\end{abstract}

\noindent\textbf{Keywords:}
hat-guessing game; infinite game; infinite combinatorics; axiom of choice; well-order

\noindent\textbf{MSC 2020:} 03E05, 91A46

\section{Introduction}

Hat-guessing games have become a rich source of problems connecting combinatorics, information theory, and set theory. A classical finite version of the problem is folklore and frequently appears in educational and recreational mathematics.

\begin{Prob} \label{original}
There are $10$ prisoners standing in a line, all facing towards the same end of the line. Each prisoner is given a hat which is either red or blue. Each prisoner can see all of those prisoners' hats who are in front of them. They guess the colour of their own hat in the order they are standing, starting with the prisoner who sees everyone else. They hear each other's guesses and can use this information to guess their own colour. The question is to determine the maximum number of correct guesses that can be guaranteed, and to give a strategy achieving it. The prisoners can fix a strategy before they get their hats, but must not communicate after that (except for hearing each other's guesses).
\end{Prob}

\begin{proof}[Solution]
We will give a strategy in which all prisoners guess correctly, except possibly for the first prisoner. The first prisoner hints with a guess of red if he sees an even number of red hats, while he hints with a guess of blue if he sees an odd number of red hats. His own guess may not be correct, but then the next prisoner can decide whether the number of red hats he sees corresponds to the first guess or not. Based on this, he can indeed guess his own hat colour. Now the next prisoners count not only the colours of seen hats, but also all the heard colours (starting from the second prisoner) and compare the parity with the hinting first guess, hence obtaining their own hat colour.
\end{proof}

The key idea is that the first prisoner uses his own guess in order to transmit one bit of global information, the parity of the visible red hats. A natural question is whether a similar strategy continues to work when the line of prisoners becomes infinite. Surprisingly, the answer is yes. The first infinite variant studied in the literature considers prisoners standing on the non-negative integers, where each prisoner still sees all subsequent prisoners but cannot hear any previous guesses. We shall call such prisoners \emph{deaf}.

\begin{Th}[Gabay, O'Connor; Hardin, Taylor \cite{HardinTaylor2013}] \label{gabay}
Assume that deaf prisoners stand on the non-negative integers facing $+\infty$, and let the colours be $\{0,1\}$. Then it is possible to guarantee that only finitely many prisoners are incorrect.
\end{Th}

\begin{proof}
We sketch a strategy guaranteeing finitely many mistakes. The prisoners can define an equivalence relation on the colour sequences with two sequences being equivalent if they differ at only finitely many places. Before getting their hats, they agree on a representative from each equivalence class. After getting the hats, each prisoner sees all but finitely many hats, thus knows the equivalence class and also its representative. The strategy for everyone is to simply guess the colour which is on themselves in the representative, therefore the guesses will form the representative itself. This means that the number of mistakes is finite as the hat sequence differs from its representative at only finitely many places.
\end{proof}

This elegant argument introduces one of the central ideas in infinite hat-guessing games: identifying colourings that differ only at finitely many places and selecting representatives of the resulting equivalence classes. The existence of such representatives relies on the axiom of choice. The following theorem further illustrates that the existence of certain hat-guessing strategies can depend on the chosen set-theoretic assumptions.

\begin{Th}[Hardin, Taylor \cite{HardinTaylor2013}]
Assume the axioms of $\textrm{ZF}$ and the Dependent Choice and the Property of Baire. Let the colours be $\{0,1\}$ and let there be countably infinitely many deaf prisoners seeing all the other hats. Then for every strategy, there is a hat sequence such that infinitely many prisoners are incorrect.
\end{Th}

Beginning with the work of Gabay and O'Connor and the work of Hardin and Taylor \cite{HardinTaylor2008}, \cite{HardinTaylor2013}, numerous simultaneous (deaf) variants have been investigated under different visibility assumptions, colour sets, and standing sets, see for example \cite{Ebert2003}, \cite{Eldredge2025}, \cite{Krzywkowski2010}. For simultaneous guessing with one-way visibility, the optimal error bounds on subsets of the real line follow from the prediction results of Hardin and Taylor, also see \cite{Hardin2013Universality} and \cite{HardinTaylor2010Minimal}.

\begin{Th}[Hardin, Taylor \cite{HardinTaylor2013}] \label{thm:deaf_one_way}
Let $S\subseteq \mathbb{R}$ be an arbitrary infinite standing set and let $C$ be a set of an arbitrary cardinality $\kappa\geq 2$. Also assume that the prisoners are deaf and each prisoner sees precisely the hats at greater positions. Then the following hold.
\begin{enumerate}
    \item If $S$ is well-ordered in its inherited order, then it can be guaranteed that only finitely many prisoners are incorrect and, for every fixed $n\in\mathbb{N}$, it cannot be guaranteed that at most $n$ prisoners are incorrect.
    \item If $S$ is not well-ordered in its inherited order, then it can be guaranteed that at most countably many prisoners are incorrect and it cannot be guaranteed that only finitely many prisoners are incorrect.
\end{enumerate}
\end{Th}

The preceding examples belong to the class of \emph{simultaneous} hat-guessing games, in which every prisoner must make their decision without access to any previous guesses. In this paper we depart from the simultaneous setting, generalising the \emph{sequential} hat-guessing games to prisoners on an arbitrary standing set $S\subseteq \mathbb{R}$. Each prisoner hears the guesses of all preceding prisoners and also sees all prisoners standing later in the ordering. This framework simultaneously encompasses several models that have appeared in the literature. It includes the sequential infinite game studied by Hardin and Taylor \cite{HardinTaylor2013}, as well as later variants considered by Glivický \cite{Glivicky2018} and Shizuma \cite{Shizuma2025}, by choosing the appropriate ordering and visibility rules. Another sequential version of the problem was investigated in \cite{Paterson2010}.

Allowing arbitrary linear orders as standing sets introduces a new phenomenon. When we have a well-ordering, the game evolves in the expected recursive fashion. Each prisoner has a well-defined turn, and the sequence of guesses is uniquely determined by the agreed strategy and the hat colours. For arbitrary linear orders, however, we found an example for the following theorem, which is counter-intuitive in our opinion. 

\begin{Th} \label{thm:example}
For the standing set $[0,1]$ and colour set $\{0,1\}$, there exists a well-defined strategy and a hat sequence, for which multiple guess sequences are valid according to that strategy on that hat sequence.
\end{Th}

Note that a similar failure of outcome uniqueness for a fixed strategy profile occurs in continuous-time repeated games, see \cite{AlosFerrerKern2015}. Consequently, one must first address the question of whether a strategy actually implies a unique outcome for each hat sequence. Despite these difficulties, our first main result shows that backwards well-ordered sets (including those which do not have a first element) retain a pleasant behaviour, familiar from the finite game.

\begin{Th} \label{thm:backwards_well_ordered}
Let $S\subseteq \mathbb{R}$ be any backwards well-ordered standing set, furthermore let the colours be $\{0,1\}$. Then it is possible to guarantee that at most $1$ prisoner is incorrect.
\end{Th}

Note that the same result for well-ordered sets was proved by Glivický \cite{Glivicky2018} for arbitrary colour sets.

This naturally leads to the following conjecture, which we believe is one of the central open problems for sequential hat-guessing games.

\begin{samepage}
\begin{Conj} \label{conj:general}
Every set $S\subseteq \mathbb{R}$ (with the induced linear order) admits a strategy with at most $1$ incorrect guess for any colour set.
\end{Conj}
\end{samepage}

Our second contribution concerns the \emph{full-visibility} model, in which every prisoner sees the hats of all other prisoners, while still guessing in the prescribed linear order. Although such a model was introduced in a paper by Glivický \cite{Glivicky2018}, there only the case of standing sets with a first element was resolved. We solve the full-visibility variant completely by proving the following general statement.

\begin{restatable}{Th}{restateFullVisibility}
\label{thm:full_visibility}
Let prisoners with full visibility stand on any set $S\subseteq \mathbb{R}$ and let the colours be chosen from a set with arbitrary cardinality. Then it is possible to guarantee that at most $1$ prisoner is incorrect.
\end{restatable}

In \cref{sec:definitions}, we formulate the sequential problem precisely and introduce the definitions and methods used throughout the paper. We also present the previously known results, stating and proving them in a framework that is consistent with our approach. In \cref{sec:original}, we prove \cref{thm:backwards_well_ordered} and provide examples that can be resolved using simpler versions of our methods. In \cref{sec:fullvisibility}, we introduce the full-visibility variant and establish the general result stated in \cref{thm:full_visibility}.

\section{Definitions and known results} \label{sec:definitions}

We begin by introducing the general version of the hat-guessing game studied throughout the paper, together with the notation used. For basic set-theoretic results and notation, we refer the reader to Jech \cite{Jech2003}. We only highlight those definitions and notations which are not standard.

\begin{Prob} \label{general}
Assume that prisoners stand on a nonempty subset $S$ of the real line. Each prisoner is given a hat which has a colour chosen from a given set $C$ of either finite or infinite cardinality $\kappa\geq 2$. Note that $C$ is also known to the prisoners before handing out the hats. Each prisoner can see those prisoners' hats who stand on greater numbers. They guess the colour of their own hat in the same order as they are standing: each prisoner hears all the guesses of those who stand on smaller numbers and can use this information to guess their own colour. We are interested in the minimal number of incorrect guesses that can be guaranteed by giving a strategy achieving it. The prisoners can fix a strategy before they get their hats, but must not communicate after that beyond hearing each other's guesses.
\end{Prob}

If $S\subseteq \mathbb{R}$, then the standing positions can be interpreted as time moments. As time flies, a prisoner standing at time $t_0$ hears all guesses made before $t_0$ and sees the hats after $t_0$. Thus, a strategy assigns a guessing function to each prisoner.

Throughout the paper, we use the term \emph{$S$-indexed sequence} for a family indexed by the standing set $S$. Such families are written in boldface. The sequence of assigned hats is denoted by $\mathbf{c}=(c_t)_{t\in S}$, and the sequence of guesses by $\mathbf{g}=(g_t)_{t\in S}$.

If the standing set $S$ has a first element, we denote the set obtained by removing that first element by $S^*$. Otherwise, we simply let $S^*=S$. Correspondingly, $\mathbf{c}^*$ and $\mathbf{g}^*$ denote the restriction of the hat sequence and guess sequence to $S^*$, respectively. In general, sequences on $S^*$ will usually be denoted with starred boldface letters, such as $\mathbf{a}^*$.

Many of our strategies depend only on the part of the hat sequence indexed by $S^*$. Therefore, throughout the paper we use the following notions for sequences on $S^*$.

Let us introduce the equivalence relation in which two colour sequences are equivalent if and only if they differ at only finitely many places. By the axiom of choice, the prisoners can choose one \emph{representative} from each equivalence class (before getting their hats). For a sequence $\mathbf{x}$ on $S^*$, let $r(\mathbf{x})$ denote the representative of the equivalence class containing $\mathbf{x}$.

If $\mathbf{x}$ and $\mathbf{y}$ differ at only finitely many places, then their \emph{Hamming distance} is denoted by $d(\mathbf{x},\mathbf{y})$, that is, the number of places where they differ. 

If $C=\{0,1\}$, let $d_2(\mathbf{x},\mathbf{y}) \in\mathbb Z_2$ denote the parity of the number of differences, that is $d_2(\mathbf{x},\mathbf{y}) = d(\mathbf{x},\mathbf{y})\pmod{2}$, and define the \emph{parity} of a sequence by $p(\mathbf{x}) = d_2(\mathbf{x},r(\mathbf{x}))$ to be the parity of the number of differences compared to its representative. Note that every representative has parity $0$. Moreover, whenever $\mathbf{x}$ and $\mathbf{y}$ differ at only finitely many places, we have $d_2(\mathbf{x},\mathbf{y}) = p(\mathbf{x})-p(\mathbf{y}) \pmod{2}$. Furthermore, observe that if a prisoner knows all but finitely many of the true hat colours, then he also knows the equivalence class, and therefore also the representative $r(\mathbf{c}^*)$ of the true hat sequence $\mathbf{c}^*$.

\begin{Def}
For a prisoner standing on $t\in S$, let $\mathbf{m}^t$ denote the \emph{mixed sequence} (or \emph{heard-seen sequence}), which is a sequence of colours on $S\setminus \{t\}$ consisting of all information available to prisoner $t$. This consists of the guesses heard before $t$ together with the hat colours seen after $t$. A mixed sequence $\mathbf{m}^t$ may be \textit{completed} to a sequence on $S$ by the guess $g_t$ of prisoner $t$, such a completion will be denoted by $\mathbf{M}^t$. We also define $r(\mathbf{m}^t)$ to be equal to $r(\mathbf{M}^t)$ for an arbitrary completion.
\end{Def}

Therefore, a \emph{strategy} for the prisoner on $t$ assigns a guess colour to each possible mixed sequence $\mathbf{m}^t$. Let us denote such a strategy by $F_t: C^{S\backslash \{t\}}\to C$. Let us denote the \emph{strategy profile} of all prisoners by $\mathcal{F}=(F_t)_{t\in S}$. For convenience, we will usually refer to the strategy profile $\mathcal{F}$ simply as the strategy.

\begin{Def}
Let a strategy $\mathcal{F}$ be fixed. We say that a guess sequence $\mathbf{g}$ is \emph{compatible} with a hat sequence $\mathbf{c}$, if for every prisoner $t\in S$, the guess $g_t$ is exactly the value prescribed by the strategy applied to the mixed sequence $\mathbf{m}^t$ determined by $\mathbf{c}$ and the previous guesses. In other words, $\mathbf{g}$ is compatible with $\mathbf{c}$ if and only if for every $t\in S$, we have $g_t=F_t(\mathbf{m}^t)$.
\end{Def}

As stated in \cref{thm:example}, unlike in the classical countable well-ordered setting, a well-defined strategy does not necessarily determine a unique compatible guess sequence. We will now prove this counter-intuitive phenomenon by giving a precise example for it.

\begin{proof}[Proof of \cref{thm:example}]
Let $S=[0,1]$ and therefore $S^*=(0,1]$, furthermore let $C=\{0,1\}$. We will try to argue with a similar strategy as in \cref{gabay}. The prisoner on $0$ will hint the parity of the true hat sequence. In the case of $S=\mathbb{N}$, this would indeed result in a strategy with at most $1$ mistake (from the fact that everyone immediately knows the correct representative), as can be seen later in \cref{thm:naturals_binary}.

At first glance, we could think that the same works for $S=[0,1]$. On the other hand, we will see that if we choose the parity function unluckily, then there may be multiple valid guess sequences for the same hat sequence.

For every $s\in S^*$, let us define the hat sequences $\mathbf{a}^s=(a^s_t)_{t\in S}$ by
\begin{equation*}
a^s_t=
\begin{cases}
    0 \textrm{ if } t\leq s,\\
    1 \textrm{ if } t > s.
\end{cases}
\end{equation*}
Furthermore, let us denote the constant $1$ sequence by $\mathbf{a}^0\equiv 1$. Notice that these sequences fall into different equivalence classes as they pairwise differ at an interval of the form $(s_1, s_2]$. Therefore we can choose each of them to be the representative of its own equivalence class. For the other equivalence classes, the representatives can be anything, let us choose them arbitrarily. Also, assume that the true hat sequence is $\mathbf{c}=\mathbf{a}^0\equiv 1$, therefore the prisoner on $0$ will hint the parity $0$.

Now observe, that $\mathbf{g}^*=\mathbf{c}^*\equiv 1$ is a valid guess sequence, as at each prisoner on $t>0$, we have $p\left((\mathbf{a}^0)^*\right)=0$, so indeed $g_t=1$. But note that everyone guessing $0$ is also a valid guess sequence, as for an arbitrary prisoner on $t>0$, we have $\mathbf{m}^t$ agreeing with $\mathbf{a}^t$ on $S\setminus \{t\}$, therefore to obtain the representative $\mathbf{a}^t$ from his point of view, he has to guess the colour $0$.
\end{proof}

The difficulty of \cref{general} in general comes from the previous observation, that even with a well-defined strategy profile $\mathcal{F}$ and a fixed hat sequence $\mathbf{c}$, there may be multiple compatible guess sequences. Because of this, throughout the paper we regard a strategy as a solution of \cref{general} only if every hat sequence admits a unique compatible guess sequence.

Before presenting the known strategies, we recall a simple observation showing that guaranteeing no incorrect guesses is impossible.

\begin{Prop} \label{non-zero}
For any nonempty standing set $S$ and any colour set $C$ with at least two elements, it cannot be guaranteed that all prisoners are correct.
\end{Prop}

\begin{proof}
Suppose for a contradiction that there exists such a strategy, for which no matter what the colour sequence $\mathbf{c}$ is, we always have a perfect guess sequence $\mathbf{g}=\mathbf{c}$. Therefore, if the colour sequence is the constant $\mathbf{c}\equiv c_1$ sequence, then everyone should guess $c_1$. On the other hand, let us change the colour of a fixed prisoner standing on $t\in S$ to a different colour $c_2$, while keeping the other colours as $c_1$. Notice that every prisoner on numbers less than $t$ should still guess $c_1$ to avoid any mistakes. But then the prisoner on $t$ heard only guesses $c_1$ and sees only colours $c_1$, so he has the same information as in the first case, so he should still guess $c_1$, which would result in a mistake.
\end{proof}

\begin{Rem}
To familiarise the reader with our notation, the proof of \cref{non-zero} can also be told in the following shorter way. Notice that for any prisoner on $t$ and hat sequence $\mathbf{c}$, changing only the colour on $t$ does not change $\mathbf{m}^t$ if we assume the correctness of each guess. As the strategy $F_t$ only depends on $\mathbf{m}^t$, we get that $g_t$ does not depend on $c_t$, giving a contradiction.
\end{Rem}

Consequently, whenever we present a strategy guaranteeing at most one incorrect guess, \cref{non-zero} immediately implies that the strategy is optimal.

We now introduce a sum operation for colours that allows us to treat finite and infinite colour sets uniformly. For this, we make use of the following well-known fact.

\begin{Lemma} \label{lem:abelian}
For an arbitrary set $A$ of (either finite or infinite) cardinality $\kappa\geq 2$, there exists an abelian group $(G,+)$ with the same cardinality $\kappa$. In the case of finite $A$, we can have $G=\mathbb{Z}_{|A|}$.
\end{Lemma}

With \cref{lem:abelian} in mind, we can always identify $C$ with its assigned abelian group. Therefore we can actually think about the colour set $C$ as an abelian group, with the colours being group elements that commute with respect to the operation $+$. Hence we can also define the difference of two colours by $c_1-c_2=c_1+(-c_2)$.

\begin{Def} \label{def:sigma_sum}
For an arbitrary nonempty set $X$, let
\[
[X]^{<\omega}
=
\{(x_1,\ldots,x_n):n\in\mathbb Z_{\ge0},\,x_i\in X\}
\]
be the set of finite ordered tuples of elements of $X$, including the empty tuple.
\end{Def}

\begin{Rem}
Throughout the paper we also use the notation $[S]^{<\omega}$ for the family of finite subsets of a standing set $S$. Thus, depending on the underlying set, the notation $[X]^{<\omega}$ refers either to the finite ordered tuples (when $X=C$) or finite subsets (when $X=S$).
\end{Rem}

\subsection{Well-ordered standing sets}

We now turn to the case where the standing set is well-ordered. Although this was settled by Glivický \cite{Glivicky2018}, who proved that at most one incorrect guess can be guaranteed for arbitrary colour sets, our framework differs substantially from theirs. So analysing this case serves to introduce our methods and provide intuition before addressing the general setting.

\begin{Th} \label{thm:finite}
Let $S$ be an arbitrary nonempty finite set, furthermore let $C$ be a set of an arbitrary cardinality $\kappa\geq 2$. Then it is possible to guarantee that at most $1$ prisoner is incorrect and it cannot be guaranteed that $0$ prisoners are incorrect.
\end{Th}

\begin{proof}
We give a strategy guaranteeing that everyone, except possibly the first prisoner, guesses correctly.

The first prisoner simply guesses the sum of all the other hat colours. The others can simply guess the difference of the hinted colour and the sum of the other prisoners' hats, which they either see or hear correctly by finite induction.
\end{proof}

\begin{Rem}
In the case of finite $S$ and infinite $C$, we could simply take a bijection between $C$ and $[C]^{<\omega}$ for the guess of the first prisoner. After this, the others would immediately get to know their colours, even without listening to the earlier guesses.
\end{Rem}

We next consider the first infinite standing set.

\begin{Th}[Hardin, Taylor \cite{HardinTaylor2013}] \label{thm:naturals_binary}
Let $S$ be the set of non-negative integers, furthermore let $C=\{0, 1\}$. Then it is possible to guarantee that at most $1$ prisoner is incorrect and it cannot be guaranteed that $0$ prisoners are incorrect.
\end{Th}

\begin{proof}
We give a strategy guaranteeing that everyone, except possibly the prisoner standing on $0$, guesses correctly.

As the prisoner on $0$ sees every other prisoner's hat, he knows the whole $\mathbf{c}^*$, from which he can calculate the representative $r(\mathbf{c}^*)$ and the parity $p(\mathbf{c}^*)$. For the strategy of the prisoner on $0$, let his guess be $g_0=p(\mathbf{c}^*)$. For the other prisoners' strategy, notice that changing only one guess also changes $p(\mathbf{g}^*)$. Therefore we can define the strategy for the prisoner on $k>0$ to be the guess which completes $\left(\mathbf{m}^k\right)^*$ to have the hinted parity $g_0=p(\mathbf{c}^*)$. Now we will prove inductively that indeed all prisoners on $k>0$ guess correctly.

First, notice that each prisoner sees all but finitely many hats, thus knows the correct $r(\mathbf{c}^*)$ and $g_0=p(\mathbf{c}^*)$. For the base case, notice that the prisoner on $1$ can indeed decide which colour completes $\left(\mathbf{m}^1\right)^*$ to the parity equalling $p(\mathbf{c}^*)$. But by definition, this is indeed $c_1$, as $\left(\mathbf{m}^1\right)^*$ and $\mathbf{c}^*$ are identical everywhere after $1$.

Now inductively, assume that everyone until $k$ guessed correctly, we want to prove that the prisoner on $k+1$ also guesses correctly. He also knows $r(\mathbf{c}^*)$ and knows that every guess before him was correct. So similarly, he can decide which colour completes $\left(\mathbf{m}^{k+1}\right)^*$ to the parity equalling $g_0=p(\mathbf{c}^*)$. Thus by induction, $\left(\mathbf{m}^{k+1}\right)^*$ and $\mathbf{c}^*$ are identical (except at $0$ and $k+1$, where one of them is not defined), so the prisoner on $k+1$ also guesses correctly, completing the induction step.
\end{proof}

The same idea extends from two colours to arbitrary colour sets.

\begin{Th}[Hardin, Taylor \cite{HardinTaylor2013}] \label{thm:naturals_kappa}
Let $S$ be the set of non-negative integers, furthermore let $C$ be a set of an arbitrary cardinality $\kappa$. Then it is possible to guarantee that at most $1$ prisoner is incorrect and it cannot be guaranteed that $0$ prisoners are incorrect.
\end{Th}

\begin{proof}
We give a strategy guaranteeing that everyone guesses correctly, except possibly for the prisoner standing on $0$.

As in \cref{thm:naturals_binary}, everyone knows the correct representative $r(\mathbf{c}^*)$. The first prisoner should guess the sum of all the differences $(\mathbf{c}^*)_k-\left(r(\mathbf{c}^*)\right)_k$ (on the places where $\mathbf{c}^*\neq r(\mathbf{c}^*)$). This is indeed well-defined, as there are only finitely many places where $\mathbf{c}^*\neq r(\mathbf{c}^*)$. Again inductively, each prisoner $t$ can simply locate the places where $\mathbf{m}^t$ does not agree with $r(\mathbf{c}^*)$, sum up the differences, and subtract this from the hinting guess $g_0$, resulting in their own difference compared to $r(\mathbf{c}^*)$.
\end{proof}

The proof of \cref{thm:naturals_binary} extends naturally to arbitrary well-ordered standing sets. Recall that every well-ordered subset of $\mathbb{R}$ is countable. Indeed, if $x\in S$ is not the last element, then it has a unique successor $x'$, so the intervals $(x,x')$ are pairwise disjoint, each containing a distinct rational number.

\begin{Th} \label{thm:wellordered_binary}
Let $S\subseteq \mathbb{R}$ be nonempty and well-ordered, furthermore let $C=\{0, 1\}$. Then it is possible to guarantee that at most $1$ prisoner is incorrect and it cannot be guaranteed that $0$ prisoners are incorrect.
\end{Th}

\begin{proof}
We give a strategy guaranteeing that everyone guesses correctly, except possibly for the first prisoner. 

The case of finite $S$ was solved in \cref{thm:finite}, so we can assume that $S$ is infinite. The strategy is the same as in \cref{thm:naturals_binary}, so the first prisoner guesses $p(\mathbf{c}^*)$, while the others guess to have $p(\mathbf{g}^*)$ being the same as the first guess $p(\mathbf{c}^*)$.

As $S$ is well-ordered, it defines an ordinal $\alpha$ such that the order of prisoners is isomorphic to $\{\beta:\beta<\alpha\}$. We prove that the strategy is well-defined and correct simultaneously by transfinite induction.

Let $\beta$ be any ordinal corresponding to a prisoner on $S^*$. By induction, we know that every prisoner on $S^*$ before $\beta$ guessed correctly. Therefore $\mathbf{m}^\beta$ and $\mathbf{c}$ are identical where they are defined. Hence from the first guess, the prisoner on $\beta$ knows the correct parity $p(\mathbf{c}^*)$, therefore he can decide which colour completes $\mathbf{m}^\beta$ to this parity, resulting in the correct guess.
\end{proof}

The transition to arbitrary colour sets is the same as in the case of transitioning from \cref{thm:naturals_binary} to \cref{thm:naturals_kappa}.

\begin{Th} \label{thm:wellordered_kappa}
Let $S\subseteq \mathbb{R}$ be nonempty and well-ordered, furthermore let $C$ be a set of an arbitrary cardinality $\kappa$. Then it is possible to guarantee that at most $1$ prisoner is incorrect and it cannot be guaranteed that $0$ prisoners are incorrect.
\end{Th}

\begin{proof}
Notice that the proof of \cref{thm:wellordered_binary} by transfinite induction works here exactly the same way, with the first prisoner guessing the sum of all the differences $(\mathbf{c}^*)_k-\left(r(\mathbf{c}^*)\right)_k$.
\end{proof}

The next theorem shows that in the case of well-ordered standing sets and the binary colour set, the strategy above using parity functions is essentially the unique optimal strategy.

\begin{Th} \label{thm:uniqueness}
Every strategy that solves \cref{thm:wellordered_binary} is essentially the same as the one we gave above. Specifically, it is true that any strategy for \cref{thm:wellordered_binary} is in the form of the first prisoner telling the parity of the number of differences for some fixed representatives.
\end{Th}

\begin{proof}
Let us fix an arbitrary colour sequence and a strategy guaranteeing at most $1$ mistake. By our assumption, we know that everyone starting from the second prisoner must guess correctly. Let us denote the guess of the first prisoner by $c\in\{0,1\}$. Switching the colour of an arbitrary but fixed non-first prisoner implies that the guess of the first prisoner should also be changed to $1-c$. Induction on the number of differences implies that for colour sequences in the same equivalence class, the guess of the first prisoner only depends on the parity of the number of differences. This determines the strategy for the whole equivalence class.

The previous reasoning works independently for each equivalence class. So for every class, the strategy is in the desired form. Finally, the choice of representatives is also independent between the classes. Note that if we choose a sequence as a representative, then changing it to another one differing in an even number of positions, the strategy remains the same.
\end{proof}

\begin{Rem}
The proof for well-ordered standing sets relies crucially on transfinite induction, and therefore does not extend to arbitrary standing sets such as $S=[0,1]$. In particular, one can no longer assume that all earlier prisoners guessed correctly. Moreover, as illustrated by \cref{thm:example}, and later by \cref{rem:equivalence}, the choice of representatives (equivalently, the parity function) may affect whether a strategy admits a unique compatible guess sequence.
\end{Rem}

\section{Backwards well-ordered standing sets} \label{sec:original}

The goal of this section is to prove \cref{thm:backwards_well_ordered} for the sequential hat-guessing game, showing that if the standing set $S$ is backwards well-ordered and $C=\{0,1\}$, then it is possible to guarantee at most one incorrect guess. To illustrate the ideas leading to the proof of the general case, we begin by proving the result for two specific sets $S$.

The sequence on $S$ (or $S^*$) that is $1$ at position $x$ and $0$ elsewhere is denoted by $\mathbbm{1}_x$.

\subsection{Two insightful subcases}

\begin{Th} \label{thm:0_1/n}
Let $S=\{0\}\bigcup \{\frac{1}{n}: n\in\mathbb{Z}_{>0}\}$, furthermore let $C=\{0, 1\}$. Then it is possible to guarantee that at most $1$ prisoner is incorrect and it cannot be guaranteed that $0$ prisoners are incorrect.
\end{Th}

\begin{proof}
We give a strategy guaranteeing that everyone, except possibly the prisoner standing on $0$, guesses correctly. The strategy is the same as in \cref{thm:naturals_binary}: the prisoner on $0$ observes the hat sequence $\mathbf{c^*}=(\ldots,c_{1/2},c_{1/1})$, and guesses $p(\mathbf{c^*})$, that is, the parity of the number of differences between the hat sequence on $S^*$ and its representative. For each $n\in \mathbb{Z}_{>0}$, the prisoner on $1/n$ observes the restriction of their mixed sequence to $S^*$, which is $\left(\mathbf{m}^{1/n}\right)^*=(\ldots, g_{1/(n+1)}, \_, c_{1/(n-1)}, \ldots, c_1)$, and guesses according to the completion $\mathbf{M}^{1/n}$ satisfying $g_0=p\left((\mathbf{M}^{1/n})^*\right)$.

It is easy to see that $\mathbf{g^*}=\mathbf{c^*}$ gives a guess sequence $\mathbf{g}$ compatible with $\mathbf{c}$ according to this strategy, since in this case, $\left(\mathbf{M}^{1/n}\right)^*$ is equal to $\mathbf{c^*}$ for all $n$, given that $g_0=p(\mathbf{c^*})$.

We show that no other guess sequences are possible. Suppose to the contrary that $\mathbf{g^*}=(\ldots,g_{1/2},g_{1/1})$ is distinct from $\mathbf{c^*}$, and take any $n$ such that $g_{1/n}\ne c_{1/n}$. Take the two completions $\mathbf{a}:=\mathbf{M}^{1/(n+1)}$ and $\mathbf{b}:=\mathbf{M}^{1/n}$. Then by the given strategy, $p(\mathbf{a^*})=p(\mathbf{b^*})=g_0$. On the other hand, observe that the sequences $\mathbf{a^*}$ and $\mathbf{b^*}$ differ only at $1/n$. In particular, both sequences have $g_{1/\ell}$ at all positions $1/\ell<1/n$ and $c_{1/\ell}$ at all positions $1/\ell>1/n$, while at position $1/n$, $\mathbf{a^*}$ has $c_{1/n}$ and $\mathbf{b^*}$ has $g_{1/n}$. So $d(\mathbf{a^*},\mathbf{b^*})=1$, contradicting the fact that the two sequences have the same parity.
\end{proof}

In this setting, it was particularly useful that $S$ had a smallest element, as this prisoner could indicate the parity of the entire hat sequence on $S^*$. We will now move to the case $S=\{\frac{1}{n}: n\in \mathbb{Z}_{>0}\}$, where this strategy does not work due to the lack of a smallest element. Instead, the finite set $D$ of positions where the hat sequence $\mathbf{c}$ differs from $r(\mathbf{c})$ will be encoded by a single prisoner $F(D)$, who will deviate from the representative of the observed sequence in order to signal $D$ to all subsequent prisoners. We will call this deviation a \textit{hint}. For the hinter to notify everyone in $D$, it is crucial that our \textit{hinter function} $F$ satisfies $F(D)<\min (D)$. Since in our proof, the prisoners on $S_{\text{even}}:=\{\frac{1}{n}: n\text{ even}\}$ will receive a different role, we also require $F$ to take values only in $S_{\text{odd}}:=\{\frac{1}{n}: n\text{ odd}\}$.

\begin{Lemma} \label{lem:hintermap1}
There exists an injective function $F: [S]^{<\omega}\to S_{\text{odd}}$ satisfying $F(D)<\min (D)$ for any nonempty $D\in {[S]}^{<\omega}$.
\end{Lemma}

\begin{proof}
Since $S$ is countably infinite, so is $[S]^{<\omega}$. Enumerate the members of $[S]^{<\omega}$ as $D_1, D_2, \ldots$ and assign inductively to each $D_i$ an element $F(D_i)\in S_{\text{odd}}$ which is smaller than $\min D_i$ and is distinct from $F(D_j)$ for all $j<i$. At each step, only finitely many elements of $S_{\text{odd}}$ are excluded as possible values of $F(D_i)$, so this assignment is possible at each step.
\end{proof}

\begin{Th} \label{thm:1/n}
Let $S=\{\frac{1}{n}: n\in\mathbb{Z}_{>0}\}$, furthermore let $C=\{0, 1\}$. Then it can be guaranteed that at most $1$ prisoner is incorrect and it cannot be guaranteed that $0$ prisoners are incorrect.
\end{Th}

\begin{proof}
Each prisoner is assigned one of the following two roles: prisoners in $S_{\text{odd}}$ are the hinters, whereas prisoners in $S_{\text{even}}$ are the alarmers. Take a function $F$ as in \cref{lem:hintermap1}. Note that $F$ is not necessarily surjective, but for simplicity, we are still going to use the terminology of hinters for every prisoner $s\in S_{\text{odd}}$. Each prisoner observes their mixed sequence $\mathbf{m}^{1/n}=(\ldots, g_{1/(n+1)}, \_, c_{1/(n-1)}, \ldots, c_1)$, and considers its representative $r(\mathbf{m}^{1/n})$. As any two prisoners' mixed sequences differ in only finitely many elements, the sequence $r(\mathbf{m}^{1/n})$ is equal for all $n$, call it $\mathbf{r}=(\ldots,r_{1/2},r_{1/1})$. Observe that $\mathbf{r}=r(\mathbf{g})$ as well, so $\mathbf{r}$ and $\mathbf{g}$ differ only at finitely many places.

\medskip

\noindent The strategy for a hinter on $1/n$ is as follows:

\begin{itemize}
    \item If there exists an earlier hinter $1/\ell$ such that $1/\ell\in \textrm{Im}(F)$ and $g_{1/\ell}\ne r_{1/\ell}$, then take the earliest such hinter, and guess $r_{1/n}+1$ if $1/n\in F^{-1}(1/\ell)$, or guess $r_{1/n}$ otherwise.

    \item Otherwise, let $B=\{1/i: 1/i>1/n \text{ and } c_{1/i}\ne r_{1/i}\}$. Guess $r_{1/n}+1$ if $F(B)=1/n$ (we call this a \textit{hint}), or guess $r_{1/n}$ otherwise.
\end{itemize}

\medskip

\noindent The strategy for an alarmer on $1/n$ is as follows:

\begin{itemize}
    \item If there exists an earlier hinter $1/\ell$ such that $1/\ell\in \textrm{Im}(F)$ and $g_{1/\ell}\ne r_{1/\ell}$, then take the earliest such hinter, and guess $r_{1/n}+1$ if $1/n\in F^{-1}(1/\ell)$, or guess $r_{1/n}$ otherwise.
    
    \item Otherwise, let $B=\{1/i: 1/i>1/n \text{ and } c_{1/i}\ne r_{1/i}\}$. Guess $r_{1/n}+1$ if $F(B)<1/n$ (we call this an \textit{alarm}), or guess $r_{1/n}$ otherwise.
\end{itemize}

\medskip

If the hat sequence $\mathbf{c}$ is given, let $\mathbf{R}=r(\mathbf{c})$ denote its representative, and let $A$ denote the set of positions $i$ such that $c_i\ne R_i$. Finally, let $h=F(A)$ be the hinter corresponding to the set $A$. We will prove that the only valid guess sequence according to this strategy is $\mathbf{g}=\mathbf{c}+\mathbbm{1}_h$.

This guess sequence is indeed compatible with $\mathbf{c}$. We have $\mathbf{r}=\mathbf{R}$ in this case, and the earliest hinter $1/\ell$ such that $g_{1/\ell}\ne r_{1/\ell}$ is indeed $h$.

\begin{itemize}
    \item Hinters $1/n<h$ will have $B=A$ so $F(B)\ne 1/n$, so they correctly guess $r_{1/n}=c_{1/n}$.
    
    \item Alarmers $1/n<h$ will have $B=A$ so $F(B)>1/n$, so they correctly guess $r_{1/n}=c_{1/n}$.
    
    \item Hinter $1/n=h$ will have $B=A$, so $F(B)=1/n$, so they guess $r_{1/n}+1=c_{1/n}+1$.
    
    \item Prisoners $1/n>h$ follow the hint of $h$ with $F^{-1}(h)=A$, guessing $r_{1/n}+1$ if $1/n\in A$ and $r_{1/n}$ otherwise, which is always $c_{1/n}$.
\end{itemize}

\medskip

Now let us prove that no other guess sequence $\mathbf{g}$ is possible. First we prove that if $\mathbf{g}$ is compatible with $\mathbf{c}$, then $\mathbf{r}=r(\mathbf{g})$ is equal to $\mathbf{R}$, i.e., $\mathbf{r}$ and $\mathbf{c}$ are equivalent. Suppose the contrary, so $\mathbf{r}$ and $\mathbf{c}$ must differ at infinitely many places. Take a prisoner $1/i$ such that $g_{1/n}=r_{1/n}$ for all $1/n\leq 1/i$, meaning that no prisoner before $1/i$ made a hint or an alarm. From the indirect assumption, we can take $1/j<1/i$ such that $r_{1/j}\ne c_{1/j}$. Take the alarmer $1/n$ immediately preceding $1/j$ (so $n\in \{j+1,j+2\}$). Prisoner $1/n$ cannot have followed any hint, and his set $B$ contains $1/j$, meaning that $F(B)<1/j$. But then also $F(B)<1/n$ (since if $1/(j+1)$ is a hinter, they did not give a hint), meaning that prisoner $1/n$ should have given an alarm, which is a contradiction.

So indeed $\mathbf{r}=\mathbf{R}$. Then for prisoners $1/n<h$, we have that $B=A$ so they do not give a hint or an alarm. Then prisoner $h$ will give a hint, and then prisoners $1/n>h$ will follow the hint of $h$, indeed giving $\mathbf{g}=\mathbf{c}+\mathbbm{1}_h$.
\end{proof}

\begin{Rem}
Observe that in the proof of \cref{thm:1/n}, no prisoners can actually give an alarm when the described strategy is carried out. Still, the alarmers play a crucial role in ruling out the unfavourable case in which the guess sequence is not equivalent to the hat sequence. 
\end{Rem}

\subsection{The general backwards well-ordered case}

\begin{Def}
We say that $S\subseteq \mathbb{R}$ is \textit{backwards well-ordered} if every nonempty subset of $S$ has a greatest element. Note that $S$ is backwards well-ordered if and only if $-S=\{-s: s\in S\}$ is well-ordered.
\end{Def}

\begin{Rem} \label{no_preceding}
Notice that in every backwards well-ordered set $S$, every non-first element $x$ has a unique preceding element. This is because we can take the nonempty subset $S_x=\{y\in S: y<x\}$. On the other hand, the converse does not hold, consider e.g. $S=\{0\}\bigcup \{\frac{1}{n}: n\in\mathbb{Z}_{>0}\}\bigcup \{2-\frac{1}{n}: n\in\mathbb{Z}_{>0}\}$.

As well-ordered subsets of $\mathbb{R}$ are countable, taking $-S$ gives that any backwards well-ordered set $S\subseteq \mathbb{R}$ is also countable.
\end{Rem}

The strategy and proof we gave for the case $S=\{0\}\cup \{\frac{1}{n}: n\in \mathbb{Z}_{>0}\}$ in \cref{thm:0_1/n} extend analogously to arbitrary backwards well-ordered sets with a first element, as follows.

\begin{Th} \label{thm:1+back-well}
Let $S\subseteq \mathbb{R}$ be backwards well-ordered with a first element, furthermore let $C=\{0, 1\}$. Then it can be guaranteed that at most $1$ prisoner is incorrect and it cannot be guaranteed that $0$ prisoners are incorrect.
\end{Th}

\begin{proof}
Let the first prisoner stand at position $t_0$. We give a strategy guaranteeing that everyone, except possibly the prisoner standing on $t_0$, guesses correctly. The strategy is the same as in \cref{thm:naturals_binary} and \cref{thm:0_1/n}, so the first prisoner observes the hat sequence $\mathbf{c}^*$ and guesses $p(\mathbf{c}^*)$, while for each $t\in S^*$, the prisoner on $t$ observes $\mathbf{m}^t$, and completes it with his guess to a sequence $\mathbf{M}^t$ satisfying $g_{t_0}=p\left((\mathbf{M}^t)^*\right)$.

Again, clearly $\mathbf{g}^*=\mathbf{c}^*$ gives a guess sequence $\mathbf{g}$ compatible with $\mathbf{c}$, since in this case, $(\mathbf{M}^t)^*$ is equal to $\mathbf{c^*}$ for all $t\in S^*$, given that $g_{t_0}=p(\mathbf{c^*})$.

To show that no other guess sequences are possible, suppose to the contrary that $g_t\ne c_t$ for some $t\in S^*$. By \cref{no_preceding}, we know that $t$ has a predecessor $s$ in the ordering of $S$. 

First, assume that $s\ne t_0$. Then $(\mathbf{M}^s)^*$ and $(\mathbf{M}^t)^*$ differ only at position $t$, where the former sequence is $c_t$ and the latter is $g_t$. Therefore the two sequences have distinct parity, so the two parities cannot both be equal to $g_{t_0}$, giving a contradiction. In the case $s=t_0$, the same argument works with the sequence $\mathbf{c}^*$ being in the place of $(\mathbf{M}^s)^*$.
\end{proof}

For backwards well-ordered sets without a first element, we proceed analogously to the case $S=\{\frac{1}{n}: n\in \mathbb{Z}_{>0}\}$ as in \cref{thm:1/n}, but we will need to make a slight modification to the strategy. The reason is that it is no longer necessarily true that any two prisoners have only finitely many other prisoners between them, allowing for the possibility that not everyone might think of the same representative.

\begin{Th} \label{thm:1/n+back-well}
Let $S\subseteq \mathbb{R}$ be nonempty and backwards well-ordered without a first element, furthermore let $C=\{0, 1\}$. Then it can be guaranteed that at most $1$ prisoner is incorrect and it cannot be guaranteed that $0$ prisoners are incorrect.
\end{Th}

\begin{proof}
The intuition of our strategy is similar to \cref{thm:1/n}, but we have to decide which prisoners are hinters ($\mathtt{h}$) and which are alarmers ($\mathtt{a}$). For this, let us define a role function $\rho$ on the elements of $S$. What we need from $\rho:S\to \{\mathtt{h}, \mathtt{a}\}$ is that if $s$ is the prisoner directly preceding $t$, then $\rho(s)\neq \rho(t)$. We can construct $\rho$ as follows. Let us define the following equivalence relation on $S$: let $s\sim t$ if and only if one of them is finitely preceded by the other, in other words, if there exist prisoners $s_1,\ldots,s_n$ such that $s_i$ directly precedes $s_{i+1}$ for all $1\le i<n$ and $\{s_1, s_n\}=\{s,t\}$. We will call the resulting equivalence classes the \textit{blocks} of $S$. Notice that each block must be isomorphic to $\{\frac{1}{n}: n\in\mathbb{Z}_{>0}\}$, since each prisoner has a direct predecessor and $S$ does not have a first element. Also notice that the blocks inherit a total order from the usual total order on $S\subseteq \mathbb{R}$, and this is also a backwards well-order.

Now we can define the role of one prisoner in each block independently, hence uniquely defining the role of all other prisoners as well by alternation. For example, we can set the last prisoner in each block to be an alarmer.

Recall that we denote the set of finite subsets of $S$ by ${[S]}^{<\omega}$. Because of the countability of $S$, similarly to \cref{lem:hintermap1}, we can take an injective function $F: {[S]}^{<\omega} \to S$ with $\rho(F(H))=\mathtt{h}$ for all $H\in {[S]}^{<\omega}$ and $F(H)<\min (H)$ for all nonempty $H\in {[S]}^{<\omega}$. Note that $F$ is not necessarily surjective, but for simplicity, we are still going to use the terminology of hinters for every prisoner $s$ with $\rho(s)=\mathtt{h}$.

Prisoners will follow almost the same strategy as in \cref{thm:1/n}, except for a slight modification for alarmers. Each prisoner $s$ observes their mixed sequence $\mathbf{m}^s$, and considers its representative $\mathbf{r}^s:=r(\mathbf{m}^s)$. The colour at position $t$ in $\mathbf{r}^s$ is denoted by $r^s_t$. The strategy for a hinter on $s$ is as follows:
\begin{itemize}
\item If there exists an earlier hinter $t$ such that $t\in \textrm{Im}(F)$ and $g_t\ne r^s_t$, then take the earliest such hinter, and guess $r^s_s+1$ if $s\in F^{-1}(t)$, or guess $r^s_s$ otherwise. (Note that such an earliest hinter exists because $\mathbf{g}$ and $\mathbf{r}^s$ can only differ at finitely many places before $s$, given that $\mathbf{m}^s$ is equal to $\mathbf{g}$ before $s$.)

\item Otherwise, let $B=\{t: t>s \text{ and } c_t\ne r^s_t\}$. (Note that $B$ is finite because $\mathbf{m}^s$ is equal to $\mathbf{c}$ after $s$.) Guess $r^s_s+1$ if $F(B)=s$ (giving a hint), or guess $r^s_s$ otherwise.
\end{itemize}
\medskip

\noindent The strategy for an alarmer on $s$ is as follows:
\begin{itemize}
\item If there exists an earlier hinter $t$ such that $t\in \textrm{Im}(F)$ and $g_t\ne r^s_t$, then take the earliest such hinter $t$. Now check whether the hint of $t$ is consistent with the hats seen by $s$, i.e., whether $F^{-1}(t)\cap \{u\in S:  u>s\}=\{u\in S: u>s \text{ and } c_u\ne r^s_u\}$. If it is inconsistent, issue an alarm, i.e., guess $r^s_s+1$. Otherwise if $s\in F^{-1}(t)$, follow the hint of $t$, so guess $r^s_s+1$, while if $s\notin F^{-1}(t)$, guess $r^s_s$.

\item If there is no such earlier hinter, let $B=\{t: t>s \text{ and } c_t\ne r^s_t\}$. Guess $r^s_s+1$ if $F(B)<s$ (issuing an alarm), or guess $r^s_s$ otherwise.
\end{itemize}
\medskip

Note that the consistency check for alarmers was added to avoid the possibility that for distinct prisoners on $s_1$, $s_2$ the sequences $\mathbf{r}^{s_1}$ and $\mathbf{r}^{s_2}$ differ. This is needed, as opposed to \cref{thm:1/n}, it is not automatically true that all prisoners observe equivalent mixed sequences.

Given a hat sequence $\mathbf{c}$ and guess sequence $\mathbf{g}$, let $\mathbf{R}=r(\mathbf{c})$ and $\mathbf{r}=r(\mathbf{g})$. As before, let $A$ denote the set of positions $i$ such that $c_i\ne R_i$, and let $h=F(A)$ be the hinter corresponding to the set $A$. We will prove that the only possible guess sequence is $\mathbf{g}=\mathbf{c}+\mathbbm{1}_h$.

First we prove that if $\mathbf{g}$ is compatible with $\mathbf{c}$, then $\mathbf{r}$ is equal to $\mathbf{R}$, i.e., $\mathbf{g}$ is equivalent to $\mathbf{c}$. Supposing the contrary, there are infinitely many prisoners $s$ such that $g_s\ne c_s$, let $E$ be the set of such prisoners. As $S$ is backwards well-ordered, we can inductively take $e_1, e_2, e_3, \ldots$ such that $e_n$ is the largest element in $E\setminus \{e_1, \ldots, e_{n-1}\}$ for every $n\in \mathbb{Z}_{>0}$. So we have that $e_1>e_2>e_3>\ldots$ is an infinite decreasing sequence in $E$.

Let $S^+$ denote the set of elements in $S$ that are greater than $e_n$ for at least one value of $n$. Then for any $s\in S^+$, there are only finitely many prisoners greater than $s$ that made an incorrect guess, therefore $\mathbf{m}^s$ is equivalent to $\mathbf{g}$. So all prisoners $s\in S^+$ observe the representative $\mathbf{r}^s=\mathbf{r}$. Take a value $i$ such that $g_t=r_t$ for all $t\in S^+$ with $t\le e_i$, meaning that no prisoner on or before $e_i$ in $S^+$ made a hint or an alarm. Note that this is possible since $\mathbf{g}$ and $\mathbf{r}$ differ in finitely many places, which can all be bounded from below by the same $e_i$.

Letting $e=e_i$, we now know that $r_e=g_e\ne c_e$. In the block of $e$, let $j$ be the last alarmer before $e$, who must still be in $S^+$. Then either $j$ directly precedes $e$, or there is a single prisoner $k$ (a hinter) between $j$ and $e$. Prisoner $j$ did not issue an alarm, and there are two cases that can explain $g_j=r_j$:\medskip

\begin{itemize}
\item If prisoner $j$ followed an earlier hinter $t\not\in S^+$, then the hint was consistent with the hats seen by $j$, so in particular, $c_e\ne r_e$ means $e\in F^{-1}(t)$. But in that case, $e$ should have followed the same hint, guessing $g_e\ne r_e$, which is a contradiction.
\item If prisoner $j$ did not detect any earlier hinter, then his set $B=\{t: t>j \text{ and } c_t\ne r_t\}$ has $F(B)>j$. But $e\in B$ and $F(B)<\min B$, so we must have $F(B)=k$ and so $k\not\in B$. But then prisoner $k$ should have given a hint, which is again a contradiction. 
\end{itemize}
\medskip

\noindent So we have shown that $\mathbf{r}=\mathbf{R}$. In that case, all $\mathbf{m}^s$ differ from $\mathbf{g}$ in only finitely many places, so $\mathbf{r}^s=\mathbf{r}$ for all $s$ as well. Then all prisoners $s<h$ have $B=A$ so they do not give a hint or an alarm. Prisoner $h$ will give a hint, and then prisoners $s>h$ follow the hint of $h$, indeed giving $\mathbf{g}=\mathbf{c}+\mathbbm{1}_h$ as the only guess sequence compatible with $\mathbf{c}$.
\end{proof}

\begin{proof}[Proof of \cref{thm:backwards_well_ordered}]
Putting together Theorems~\ref{thm:1+back-well} and \ref{thm:1/n+back-well} finishes the proof of \cref{thm:backwards_well_ordered} in all cases.
\end{proof}

\subsection{Further remarks}

\begin{Rem} \label{rem:colours}
In \cref{sec:fullvisibility}, we will solve the version with full visibility for any standing set $S$ and any colour set $C$. The method will be similar to the one used in this section. On the other hand, we will use multiple and modified functions $F$ for assigning the hinting prisoners, to deal with the arbitrary set of colours. Using the same trick, we can extend all of the above results to any colour set $C$.
\end{Rem}

\begin{Rem} \label{rem:equivalence}
Instead of talking about abstract strategies, it might be more practical to talk about purely set-theoretic questions. We can actually give such a set-theoretic condition for the case of $S=[0,1]$ and $C=\{0,1\}$, which holds if and only if it can be guaranteed that at most $1$ prisoner is incorrect.

For the compatibility of the guess sequence in which every prisoner on $(0,1]$ is correct, the prisoner on $0$ must guess the parity of $\mathbf{c}^*$, similarly to \cref{thm:uniqueness}. Therefore, the only flexibility in the strategy is the choice of the parity function.

For colour sequences $\mathbf{a}^*$ and $\mathbf{b}^*$ on $S^*=(0,1]$, let us define their combination at $t$ to be the following sequence $\langle{\mathbf{a}^*,\mathbf{b}^*\rangle}^t$ on $(0,1]$:
\begin{equation*}
\left({\langle \mathbf{a}^*,\mathbf{b}^*\rangle}^t\right)_s=
\begin{cases}
    (\mathbf{a}^*)_s
    & \textrm{ if } 0<s\leq t,\\
    (\mathbf{b}^*)_s
    & \textrm{ if } t<s\leq 1.
\end{cases}
\end{equation*}

Let us fix the hat sequence to be $\mathbf{b}$. In order to guarantee that no guess sequence $\mathbf{a}$ with $\mathbf{a}^*\neq \mathbf{b}^*$ can be compatible, there must exist a prisoner $t$ who guessed according to the wrong parity. Precisely, there exists a strategy guaranteeing at most $1$ incorrect guess, if and only if, there exists a parity function $P$ on the colour sequences such that for any hat sequence $\mathbf{c}$ and guess sequence with $\mathbf{g}^*\neq \mathbf{c}^*$, there exists a $t\in (0,1]$ for which $P\left({\langle \mathbf{g}^*,\mathbf{c}^* \rangle}^t\right)\neq P(\mathbf{c}^*)$.
\end{Rem}

\section{The full visibility version} \label{sec:fullvisibility}

The most well-known theorem from the literature with full visibility is probably originating from Lenstra. This states that in the case of full visibility with simultaneous guessing, the prisoners can guarantee that either all of them are correct or all of them are incorrect.

\begin{Th}[Lenstra; Hardin, Taylor \cite{HardinTaylor2013}]
Let $C=\{0,1\}$ and let there be deaf prisoners of arbitrary cardinality with full visibility. We state that it is possible to give a strategy such that either all prisoners are correct, or all are incorrect.
\end{Th}

In this section, we will cover a modified version of the problem. We will keep the sequential guessing from the earlier sections, the modification is that the prisoners have full visibility, meaning that they can see all the other hats. Only the order of the guesses is based on the standing set $S$. Although such a model was introduced in a paper by Glivický \cite{Glivicky2018}, there only the case of standing sets with a first element was resolved.

The idea to include full visibility in the hearing version in our framework originates from Tamás Keleti. He also contributed to the proofs of these theorems. Now, we will give a complete solution to this setting.

\subsection{Solving the full visibility version}

Note that because of the full visibility, each prisoner knows the correct representative. So if each prisoner would guess the colour from the representative, that would result in finitely many incorrect guesses. Actually, for every standing set $S\subseteq \mathbb{R}$ and any colour set $C$, we will see that it can also be guaranteed to have only at most $1$ incorrect guess, as claimed in \cref{thm:full_visibility}.

Note that for any $S\neq \emptyset$, similarly to \cref{non-zero}, we cannot guarantee $0$ incorrect guesses. We prove \cref{thm:full_visibility} by proving it for different kinds of standing sets and colour sets.

\begin{Th} \label{thm:finite_fullvis}
Let $S\subseteq \mathbb{R}$ be an arbitrary nonempty finite standing set and let $C$ be a set of an arbitrary cardinality $\kappa\geq 2$. Also assume that the prisoners have full visibility. Then it can be guaranteed that at most $1$ prisoner is incorrect and it cannot be guaranteed that $0$ prisoners are incorrect.
\end{Th}

\begin{proof}
The proof of \cref{thm:finite} also works here, as now the prisoners have more vision and therefore more information.
\end{proof}

The following result can also be obtained from the algebraic construction of Glivický \cite{Glivicky2018}. We include a short proof here, as our general framework and methods differ from theirs.

\begin{Th} \label{thm:infinite_fullvis_1}
Let $S\subseteq \mathbb{R}$ be an arbitrary infinite standing set with a first element and let $C$ be a set of an arbitrary cardinality $\kappa\geq 2$. Also assume that the prisoners have full visibility. Then it can be guaranteed that at most $1$ prisoner is incorrect and it cannot be guaranteed that $0$ prisoners are incorrect.
\end{Th}

\begin{proof}
Note that as the prisoners have full visibility, they instantly know the true equivalence class and the representative of the hat sequence. Similarly to \cref{thm:naturals_kappa}, the first prisoner can simply guess the sum of the differences compared to the representative. After this first guess, all the other prisoners instantly get to know the colour of their hats, as they can subtract the seen differences from the hint, resulting in their own difference.
\end{proof}

In the case of $S$ without a first element, we will assign hinting prisoners to the finite subsets of $S$, similarly to \cref{thm:1/n+back-well}. Unfortunately, because of the arbitrary colour set, the strategy gets complicated if an assigned hinting prisoner would have to guess the same colour as the one from the representative, because then the later prisoners would not notice the hint. Both in the case of finite and infinite $C$, we will deal with this by having two separate hinting functions: one for the general case when the sum of differences is not the zero element, and a separate one for the zero-sum case. Furthermore, such functions do not necessarily exist in general for the whole standing set, for example for $S=\{\frac{1}{n}: n\in\mathbb{Z}_{>0}\}\bigcup [1,2]$. Because of this, we will only look for them on an initial segment of $S$. For the existence of such appropriate functions, we prove \cref{lem:hintermap2}, in a manner similar to \cref{lem:hintermap1}.

\begin{Lemma} \label{lem:hintermap2}
For any infinite $S\subseteq \mathbb{R}$ without a first element, there exists a $t\in \mathbb{R}$ such that for the slice $S_t=\{s\in S: s<t\}$, the following holds. There are injective functions $F_0,F_1: {[S_t]}^{<\omega} \to S_t$ with disjoint image (so $F_0(D_0)\neq F_1(D_1)$ for any $D_0,D_1\in {[S_t]}^{<\omega}$) such that $F_0(D)<\min (D)$ and $F_1(D)<\min (D)$ for any nonempty $D\in {[S_t]}^{<\omega}$. Furthermore, we can also assume that $T:=S\setminus S_t$ is either empty or infinite.
\end{Lemma}

\begin{proof}
Let $m=\inf(S)$, by definition $m\notin S$. Note that it might happen that $m=-\infty$. By the well-ordering of cardinals, there is an $\alpha = \min_{s\in S} \left|S \bigcap (m,s)\right|$, minimised at some $t>m$, and let $S_t=S\bigcap (m,t)$. If $S\setminus S_t$ is finite, then $S$ is bounded from above, so we may replace $t$ by a real number greater than every element of $S$ and let $S_t=S$. As $\alpha$ cannot be finite, it has to be infinite, so $\left|{[S_t]}^{<\omega}\right|=\left|S_t\right|$. From the definition of $\alpha$ and $S_t$, each nonempty subslice $S_u\subseteq S_t$ also has $|S_u|=|S_t|$. Now by well-ordering the elements of ${[S_t]}^{<\omega}$, we can indeed define the injective functions $F_0,F_1$ simultaneously recursively.
\end{proof}

\begin{Th} \label{thm:infinite_fullvis_2}
Let $S\subseteq \mathbb{R}$ be any infinite standing set without a first element and let $C$ be a set of an arbitrary cardinality $\kappa\geq 2$. Also assume that the prisoners have full visibility. Then it can be guaranteed that at most $1$ prisoner is incorrect and it cannot be guaranteed that $0$ prisoners are incorrect.
\end{Th}

\begin{proof}
First use \cref{lem:hintermap2} to get a $t\in \mathbb{R}$ and $S_t=\{s\in S: s<t\}$, with injective functions $F_0,F_1: {[S_t]}^{<\omega} \to S_t$ with disjoint image such that $F_0(D),F_1(D)<\min (D)$ for any nonempty $D\in {[S_t]}^{<\omega}$. Furthermore, we also got that $T=S\setminus S_t$ is either empty or infinite. Before getting the hats, for each colour $c$, fix a colour $n(c)\neq c$.

\textit{Case 1: $T$ is empty.} In this case we have $S=S_t$. By definition, there are finitely many places at which the representative of $S_t$ differs from the true hat sequence on $S_t$, let this finite set be $D$. Let the sum of the differences compared to the representative be $c\in C$. If $c\neq 0$, then let $F_1(D)$ be the hinting prisoner, let his colour from the representative be $c'$. Then this hinting prisoner should guess the colour $c+c'$. If $c=0$, then let $F_0(D)$ be the hinting prisoner, let his colour from the representative be $c'$. Then this hinting prisoner should guess $n(c')$.

\textit{Case 2: $T$ is infinite.} In this case, by definition, there are finitely many places at which the restriction of the representative on $S_t$ differs from the true hat sequence on $S_t$, let this finite set be $D\subseteq S_t$. Also by definition, there are finitely many places at which the restriction of the representative on $T$ differs from the true hat sequence on $T$, let this finite set be $D_T\subseteq T$. Let $c$ be the sum of all differences compared to the representative, so the sum of differences compared to the representative on $D$ and the sum of differences compared to the representative on $D_T$. If $c\neq 0$, then let $F_1(D)$ be the hinting prisoner, let his colour from the representative be $c'$. Then this hinting prisoner should guess the colour $c+c'$. If $c=0$, then let $F_0(D)$ be the hinting prisoner, let his colour from the representative be $c'$. Then this hinting prisoner should guess $n(c')$.

Note that prior to the guessing, the prisoners do not know who is the hinting prisoner. Furthermore, as the prisoners have full visibility, they know whether the earlier ones were correct or not. If a prisoner has already heard more than $1$ incorrect guess, then he should automatically guess the colour from the representative. This addition serves to guarantee that there is no other valid guess sequence beyond the correct one.

Otherwise, if a hint has already been heard, the prisoner should use the earliest such hint, meaning the earliest guess different from the representative made by an earlier prisoner in $\textrm{Im}(F_0)\cup\textrm{Im}(F_1)$. (Note that only finitely many earlier guesses can differ from the representative, since at most $1$ earlier guess is incorrect and the true colours differ from the representative at only finitely many places. Thus the earliest such hint exists. Also note that in reality, in every valid guess sequence there will be exactly $1$ hint, but this is needed to completely well-define the strategy.) If the prisoner's place belongs to the hinted set $D$ or $T$, then as he can see all the other true colours and knows the hinted set and the sum of differences, he can calculate his own true colour, otherwise he should guess the colour from the representative.

If no earlier hint has been heard, a prisoner on $S_t$ determines whether he should hint by applying the relevant case with his own hat colour taken to be its colour from the representative on $S_t$. If this calculation designates him as the hinting prisoner, he should guess according to the above cases, otherwise he should guess the colour from the representative. To have the strategy well-defined (even though the following would never occur in reality), the prisoners on $T$ who have heard no earlier hint should guess the colour from the representative.

Formally, let $\mathbf{r}=(r_t)_{t\in S}$ be obtained by combining the representatives on $S_t$ and (when nonempty) $T$. Every prisoner knows $\mathbf{r}$ by full visibility, and only finitely many differences $c_t-r_t$ are nonzero. The strategy for a prisoner on $s\in S$ is as follows:\medskip

\begin{itemize}
    \item If $\left|\{u\in S:u<s\text{ and }g_u\ne c_u\}\right|>1$, guess $r_s$.

    \item Otherwise, if there exists an earlier prisoner $h<s$ such that $h\in\textrm{Im}(F_0)\cup\textrm{Im}(F_1)$ and $g_h\ne r_h$, take the earliest such prisoner. As $h=F_i(D)$, where $i\in\{0,1\}$ and $D\in[S_t]^{<\omega}$ are uniquely determined, let $c=0$ if $i=0$, while $c=g_h-r_h$ if $i=1$. If $s\in D\cup T$, guess $r_s+c-\sum_{u\in S\setminus\{s\}}(c_u-r_u)$, otherwise guess $r_s$.

    \item Otherwise if $s\in S_t$, let $B=\{u\in S_t\setminus\{s\}:c_u\ne r_u\}$ and $c=\sum_{u\in S\setminus\{s\}}(c_u-r_u)$. If $c\neq 0$ and $F_1(B)=s$, guess $r_s+c$, while if $c=0$ and $F_0(B)=s$, guess $n(r_s)$, otherwise guess $r_s$.

    \item Otherwise ($s\in T$ and) guess $r_s$.
\end{itemize}
\medskip

\noindent We now prove the correctness of this strategy. It is easy to see that the guess sequence in which the hinting prisoner gives the prescribed hint and everyone else guesses correctly is indeed compatible. Now consider any compatible guess sequence. Among any three prisoners who guess incorrectly, the last one has heard at least two incorrect guesses, so he must guess the colour from the representative. Therefore, there can be at most two incorrect guesses differing from the representative, so altogether there are only finitely many incorrect guesses. So if there were two distinct compatible guess sequences for the same hat sequence, they would also differ at only finitely many places. Take the first prisoner whose guess differ. This prisoner sees the same hat colours and hears the same earlier guesses in both sequences, so the strategy prescribes the same guess, giving a contradiction. Therefore the only compatible guess sequence is the desired one.
\end{proof}

\begin{Rem}
Notice that in the proofs of Theorems~\ref{thm:1/n}, \ref{thm:1/n+back-well}, and \ref{thm:infinite_fullvis_2}, we gave strategies which guarantee that exactly $1$ prisoner is incorrect.
\end{Rem}

\begin{proof}[Proof of \cref{thm:full_visibility}]
Putting together Theorems~\ref{thm:finite_fullvis}, \ref{thm:infinite_fullvis_1}, and
\ref{thm:infinite_fullvis_2}, indeed proves the statement for any standing set $S$ and colour set $C$.
\end{proof}

\section*{Acknowledgements}

We would like to thank Péter Mazug, with whom we first explored these questions and whose many insightful discussions and enjoyable brainstorming sessions contributed greatly to the development of this work. We are grateful to Tamás Keleti for suggesting the full-visibility variant and for contributing to the development of the corresponding proofs.

\section*{AI declaration}

The results and underlying proof methods presented in this manuscript were developed by the authors. GPT-6 Astra was used to assist with language polishing, grammar and checking proofs. A separate forthcoming manuscript by the same authors will present a proof generated by GPT-6 Astra that completely resolves the one-way visibility problem for all standing sets $S\subseteq \mathbb{R}$, proving that at most $1$ incorrect guess can be achieved. The authors take full responsibility for the mathematical correctness and final content of this manuscript.

\printbibliography[heading=bibintoc, title={References}]

\end{document}